\documentclass[12pt]{article}

\catcode`\@=11

\thicklines
\newskip\Einheit \Einheit=.6cm
\newcount\xcoord \newcount\ycoord
\newdimen\xdim \newdimen\ydim \newdimen\PfadD@cke \newdimen\Pfadd@cke
\def\PfadDicke#1{\PfadD@cke#1 \divide\PfadD@cke by2 
\Pfadd@cke\PfadD@cke \multiply\PfadD@cke by2}
\long\def\LOOP#1\REPEAT{\def\BODY{#1}\ITERATE}
\def\ITERATE{\BODY \let\next\ITERATE \else\let\next\relax\fi \next}
\let\REPEAT=\fi
\def\Punkt{\hbox{\raise-2pt\hbox to0pt{\hss\scriptsize$\bullet$\hss}}}

\def\DuennPunkt(#1,#2){\unskip
  \raise#2 \Einheit\hbox to0pt{\hskip#1 \Einheit
          \raise-1.5pt\hbox to0pt{\hss\tiny$\bullet$\hss}\hss}}
		  
\def\NormalPunkt(#1,#2){\unskip
  \raise#2 \Einheit\hbox to0pt{\hskip#1 \Einheit
          \raise-3pt\hbox to0pt{\hss\large$\bullet$\hss}\hss}}
\def\DickPunkt(#1,#2){\unskip
  \raise#2 \Einheit\hbox to0pt{\hskip#1 \Einheit
          \raise-4pt\hbox to0pt{\hss\Large$\bullet$\hss}\hss}}
\def\Kreis(#1,#2){\unskip
  \raise#2 \Einheit\hbox to0pt{\hskip#1 \Einheit
          \raise-4pt\hbox to0pt{\hss\Large$\circ$\hss}\hss}}
\def\Diagonale(#1,#2)#3{\unskip\leavevmode
  \xcoord#1\relax \ycoord#2\relax
      \raise\ycoord \Einheit\hbox to0pt{\hskip\xcoord \Einheit
         \unitlength\Einheit
         \line(1,1){#3}\hss}}
\def\AntiDiagonale(#1,#2)#3{\unskip\leavevmode
  \xcoord#1\relax \ycoord#2\relax \advance\xcoord by -0.05\relax
      \raise\ycoord \Einheit\hbox to0pt{\hskip\xcoord \Einheit
         \unitlength\Einheit
         \line(1,-1){#3}\hss}}
\def\Pfad(#1,#2),#3\endPfad{\unskip\leavevmode
  \xcoord#1 \ycoord#2 \thicklines\ZeichnePfad#3\endPfad\thinlines}
\def\ZeichnePfad#1{\ifx#1\endPfad\let\next\relax
  \else\let\next\ZeichnePfad
    \ifnum#1=1
      \raise\ycoord \Einheit\hbox to0pt{\hskip\xcoord \Einheit
         \vrule height\Pfadd@cke width1 \Einheit depth\Pfadd@cke\hss}%
      \advance\xcoord by 1
     \else\ifnum#1=2
      \raise\ycoord \Einheit\hbox to0pt{\hskip\xcoord \Einheit
         \unitlength\Einheit
         \line(0,1){1}\hss}
      \advance\xcoord by 0
      \advance\ycoord by 1
 \else\ifnum#1=3
      \raise\ycoord \Einheit\hbox to0pt{\hskip\xcoord \Einheit
         \unitlength\Einheit
         \line(1,1){1}\hss}
      \advance\xcoord by 1
      \advance\ycoord by 1
    \else\ifnum#1=4
      \raise\ycoord \Einheit\hbox to0pt{\hskip\xcoord \Einheit
         \unitlength\Einheit
         \line(1,-1){1}\hss}
      \advance\xcoord by 1
      \advance\ycoord by -1
   \else\ifnum#1=5
      \raise\ycoord \Einheit\hbox to0pt{\hskip\xcoord \Einheit
         \unitlength\Einheit
         \line(2,1){2}\hss}
      \advance\xcoord by 2
      \advance\ycoord by 1
	  \else\ifnum#1=6
      \raise\ycoord \Einheit\hbox to0pt{\hskip\xcoord \Einheit
         \unitlength\Einheit
         \line(2,-1){2}\hss}
      \advance\xcoord by 2
      \advance\ycoord by -1
	  \else\ifnum#1=7
      \raise\ycoord \Einheit\hbox to0pt{\hskip\xcoord \Einheit
         \unitlength\Einheit
         \line(3,1){3}\hss}
      \advance\xcoord by 3
      \advance\ycoord by 1
	  \else\ifnum#1=8
      \raise\ycoord \Einheit\hbox to0pt{\hskip\xcoord \Einheit
         \unitlength\Einheit
         \line(3,-1){3}\hss}
      \advance\xcoord by 3
      \advance\ycoord by -1
    \fi\fi\fi\fi\fi\fi\fi\fi
  \fi\next}
\def\hSSchritt{\leavevmode\raise-.4pt\hbox 
to0pt{\hss.\hss}\hskip.2\Einheit
  \raise-.4pt\hbox to0pt{\hss.\hss}\hskip.2\Einheit
  \raise-.4pt\hbox to0pt{\hss.\hss}\hskip.2\Einheit
  \raise-.4pt\hbox to0pt{\hss.\hss}\hskip.2\Einheit
  \raise-.4pt\hbox to0pt{\hss.\hss}\hskip.2\Einheit}
\def\vSSchritt{\vbox{\baselineskip.2\Einheit\lineskiplimit0pt
\hbox{.}\hbox{.}\hbox{.}\hbox{.}\hbox{.}}}
\def\DSSchritt{\leavevmode\raise-.4pt\hbox to0pt{%
  \hbox to0pt{\hss.\hss}\hskip.2\Einheit
  \raise.2\Einheit\hbox to0pt{\hss.\hss}\hskip.2\Einheit
  \raise.4\Einheit\hbox to0pt{\hss.\hss}\hskip.2\Einheit
  \raise.6\Einheit\hbox to0pt{\hss.\hss}\hskip.2\Einheit
  \raise.8\Einheit\hbox to0pt{\hss.\hss}\hss}}
\def\dSSchritt{\leavevmode\raise-.4pt\hbox to0pt{%
  \hbox to0pt{\hss.\hss}\hskip.2\Einheit
  \raise-.2\Einheit\hbox to0pt{\hss.\hss}\hskip.2\Einheit
  \raise-.4\Einheit\hbox to0pt{\hss.\hss}\hskip.2\Einheit
  \raise-.6\Einheit\hbox to0pt{\hss.\hss}\hskip.2\Einheit
  \raise-.8\Einheit\hbox to0pt{\hss.\hss}\hss}}
\def\SPfad(#1,#2),#3\endSPfad{\unskip\leavevmode
  \xcoord#1 \ycoord#2 \ZeichneSPfad#3\endSPfad}
\def\ZeichneSPfad#1{\ifx#1\endSPfad\let\next\relax
  \else\let\next\ZeichneSPfad
    \ifnum#1=1
      \raise\ycoord \Einheit\hbox to0pt{\hskip\xcoord \Einheit
         \hSSchritt\hss}%
      \advance\xcoord by 1
    \else\ifnum#1=2
      \raise\ycoord \Einheit\hbox to0pt{\hskip\xcoord \Einheit
        \hbox{\hskip-2pt \vSSchritt}\hss}%
      \advance\ycoord by 1
    \else\ifnum#1=3
      \raise\ycoord \Einheit\hbox to0pt{\hskip\xcoord \Einheit
         \DSSchritt\hss}
      \advance\xcoord by 1
      \advance\ycoord by 1
    \else\ifnum#1=4
      \raise\ycoord \Einheit\hbox to0pt{\hskip\xcoord \Einheit
         \dSSchritt\hss}
      \advance\xcoord by 1
      \advance\ycoord by -1
    \fi\fi\fi\fi
  \fi\next}
\def\Koordinatenachsen(#1,#2){\unskip
 \hbox to0pt{\hskip-.5pt\vrule height#2 \Einheit width.5pt depth1 
\Einheit}%
 \hbox to0pt{\hskip-1 \Einheit \xcoord#1 \advance\xcoord by1
    \vrule height0.25pt width\xcoord \Einheit depth0.25pt\hss}}
\def\Koordinatenachsen(#1,#2)(#3,#4){\unskip
 \hbox to0pt{\hskip-.5pt \ycoord-#4 \advance\ycoord by1
    \vrule height#2 \Einheit width.5pt depth\ycoord \Einheit}%
 \hbox to0pt{\hskip-1 \Einheit \hskip#3\Einheit 
    \xcoord#1 \advance\xcoord by1 \advance\xcoord by-#3 
    \vrule height0.25pt width\xcoord \Einheit depth0.25pt\hss}}
\def\Gitter(#1,#2){\unskip \xcoord0 \ycoord0 \leavevmode
  \LOOP\ifnum\ycoord<#2
    \loop\ifnum\xcoord<#1
      \raise\ycoord \Einheit\hbox to0pt{\hskip\xcoord 
\Einheit\Punkt\hss}%
      \advance\xcoord by1
    \repeat
    \xcoord0
    \advance\ycoord by1
  \REPEAT}
\def\Gitter(#1,#2)(#3,#4){\unskip \xcoord#3 \ycoord#4 \leavevmode
  \LOOP\ifnum\ycoord<#2
    \loop\ifnum\xcoord<#1
      \raise\ycoord \Einheit\hbox to0pt{\hskip\xcoord 
\Einheit\Punkt\hss}%
      \advance\xcoord by1
    \repeat
    \xcoord#3
    \advance\ycoord by1
  \REPEAT}
\def\Label#1#2(#3,#4){\unskip \xdim#3 \Einheit \ydim#4 \Einheit
  \def\lo{\advance\xdim by-.5 \Einheit \advance\ydim by.5 \Einheit}%
  \def\llo{\advance\xdim by-.25cm \advance\ydim by.5 \Einheit}%
  \def\loo{\advance\xdim by-.5 \Einheit \advance\ydim by.25cm}%
  \def\o{\advance\ydim by.25cm}%
  \def\ro{\advance\xdim by.5 \Einheit \advance\ydim by.5 \Einheit}%
  \def\rro{\advance\xdim by.25cm \advance\ydim by.5 \Einheit}%
  \def\roo{\advance\xdim by.5 \Einheit \advance\ydim by.25cm}%
  \def\l{\advance\xdim by-.30cm}%
  \def\r{\advance\xdim by.30cm}%
  \def\lu{\advance\xdim by-.5 \Einheit \advance\ydim by-.6 \Einheit}%
  \def\llu{\advance\xdim by-.25cm \advance\ydim by-.6 \Einheit}%
  \def\luu{\advance\xdim by-.5 \Einheit \advance\ydim by-.30cm}%
  \def\u{\advance\ydim by-.30cm}%
  \def\ru{\advance\xdim by.5 \Einheit \advance\ydim by-.6 \Einheit}%
  \def\rru{\advance\xdim by.25cm \advance\ydim by-.6 \Einheit}%
  \def\ruu{\advance\xdim by.5 \Einheit \advance\ydim by-.30cm}%
  #1\raise\ydim\hbox to0pt{\hskip\xdim
     \vbox to0pt{\vss\hbox to0pt{\hss$#2$\hss}\vss}\hss}%
}
\catcode`\@=12

\usepackage{amsmath,amsthm}
\usepackage{amssymb}
 
\usepackage{color}
\usepackage{xspace}
\usepackage[colorlinks=true,
linkcolor=green,
filecolor=brown,
citecolor=green]{hyperref}

\def\blue{\textcolor{blue} }

\def\white{\textcolor{white} }

\def\gl{ground level\xspace}

\def\mbf#1{\mathchoice{\hbox{\boldmath $\displaystyle #1$}}
        {\hbox{\boldmath $\textstyle #1$}}
        {\hbox{\boldmath $\scriptstyle #1$}}
        {\hbox{\boldmath $\scriptscriptstyle #1$}}} 

\begin{document}
\newtheorem{theorem}{Theorem}
\newtheorem{defn}[theorem]{Definition}
\newtheorem{lemma}[theorem]{Lemma}
\newtheorem{prop}[theorem]{Proposition}
\newtheorem{cor}[theorem]{Corollary}
 
\begin{center}
{\Large
Another bijection for 021-avoiding ascent sequences      \\ 
}

\vspace{5mm}
David Callan  \\
\noindent {\small Dept. of Statistics, 
University of Wisconsin-Madison,  Madison, WI \ 53706}  \\
{\bf callan@stat.wisc.edu} 

February 23, 2014
\end{center}

\begin{abstract}
Chen and collaborators give a recursively defined bijection from 021-avoiding 
ascent sequences to 021-avoiding (aka 132-avoiding) permutations. Here we 
give an algorithmic bijection from 021-avoiding ascent sequences to 
Dyck paths. 
Our bijection does not appear to be closely related to the Chen bijection but, 
like the Chen bijection, it preserves several interesting statistics. 

\end{abstract}

\vspace{5mm}

\section{Introduction} 

\vspace*{-3mm}

Several recent papers treat pattern avoidance in ascent 
sequences \cite{chen, duncan, some, yan14}. 
A striking result \cite{duncan} is that 021-avoiding ascent sequences 
are counted by the Catalan numbers. William Chen and his collaborators 
\cite{chen}  give an elegant recursively defined bijection from 021-avoiding 
ascent sequences to 021-avoiding (aka 132-avoiding) permutations. In this paper 
we give a bijection from 021-avoiding ascent sequences to Dyck paths   
that builds up the Dyck path iteratively. 
In Section \ref{two} we recall the relevant definitions for ascent sequences 
and  terminology for Dyck paths. Section \ref{three} presents the bijection and 
Section \ref{four} its inverse. Lastly, Section \ref{five} mentions some 
statistics preserved by the bijection.

\section{Ascent sequences and Dyck paths} \label{two}

\vspace*{-3mm}

An \emph{ascent} in a sequence of integers is a pair of consecutive entries with 
the first smaller than the second. An \emph{ascent sequence} is a 
sequence $(u_1,u_2,\dots,u_n)$ of nonnegative integers 
such that $u_0=0$ and $u_i \le 1 + $ \# ascents in $(u_1,u_2,\dots,u_{i-1})$ 
for $i \ge 2$. Due to the initial 0, it is clear that 021-avoiding ascent 
sequences can be characterized as ascent sequences in which the nonzero entries are weakly increasing.

A Dyck path is a lattice path of upsteps $U=(1,1)$ and 
downsteps $D=(1,-1)$, the same number of each, that stays weakly 
above the horizontal line, called \emph{\gl,} that joins its initial 
and terminal points (vertices). We have the usual notions of size (number of upsteps), ascent (maximal sequence of contiguous upsteps), descent, peak ($UD$), peak vertex 
(the vertex between the $U$ and $D$), valley ($DU$), and valley vertex. 
An ascent is \emph{short} if it has length 1, otherwise \emph{long}.

\Einheit=0.6cm
\[
\hspace*{-10mm}
\Label\o{\rightarrow}(-5.2,3)
\Label\u{\uparrow}(-2,1.9)
\Label\l{ \textrm{{\footnotesize peak upstep}}}(-6.9,3.5)
\Label\u{ \textrm{{\footnotesize valley}}}(-2,1)
\Label\u{ \textrm{{\footnotesize vertex}}}(-2,0.4)
\Label\u{\uparrow}(.4,1.5)
\Label\u{ \textrm{{\footnotesize return}}}(.4,.8)
\Label\u{ \textrm{{\footnotesize downstep}}}(.4,0.2)
\SPfad(-7,1),1111\endSPfad
\SPfad(1,1),11111111\endSPfad
\SPfad(-1,1),1\endSPfad
\Pfad(-7,1),3334434433433444\endPfad
\DuennPunkt(-7,1)
\DuennPunkt(-6,2)
\DuennPunkt(-5,3)
\DuennPunkt(-4,4)
\DuennPunkt(-3,3)
\DuennPunkt(-2,2)
\DuennPunkt(-1,3) 
\DuennPunkt(0,2) 
\DuennPunkt(1,1)
\DuennPunkt(2,2)
\DuennPunkt(3,3)
\DuennPunkt(4,2)
\DuennPunkt(5,3)
\DuennPunkt(6,4)
\DuennPunkt(7,3)
\DuennPunkt(8,2)
\DuennPunkt(9,1)
\Label\u{ \textrm{{\footnotesize key downsteps}}}(9.5,4.5)
\Label\u{ \textrm{{\footnotesize (in blue)}}}(9.5,3.8)
\Label\u{\swarrow}(9.0,3.1)
\Label\u{\uparrow}(5,1)
\Label\u{ \textrm{{\footnotesize \gl}}}(5,.2)
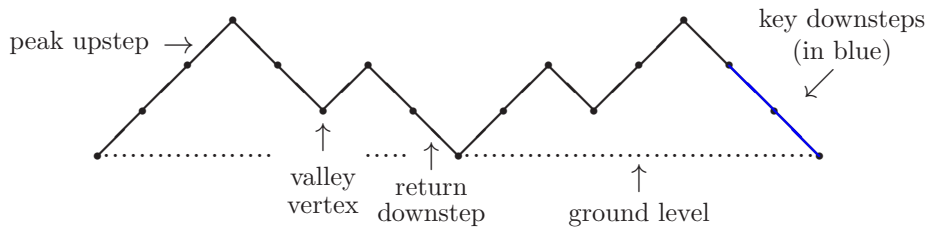
\Label\o{ \textrm{\small  A non-elevated Dyck path of size 7 with  2 key downsteps, long last ascent}}(1,-2.3)
\blue{
\Pfad(7,3),44\endPfad }
\]
  
\vspace*{2mm}

The height of a vertex in a Dyck path is its vertical height above \gl. 
A \emph{return} downstep is one that returns the path to \gl. An \emph{elevated} 
Dyck path is one with exactly one return (necessarily at the end). 
The \emph{degree of elevation} of a Dyck path is the height of its lowest valley vertex (undefined for \emph{pyramid} Dyck paths---$U^n D^n$---which have no valleys). Thus the degree of elevation is 0 precisely for non-elevated Dyck paths.
Upsteps and downsteps come in matching pairs: travel due east 
from an upstep to the first downstep encountered. More precisely, 
$D_{0}$ is the matching downstep for upstep $U_{0}$ if $D_{0}$ 
terminates the shortest Dyck subpath that starts with $U_{0}$. 
It is convenient to define a \emph{key downstep} in a Dyck path 
to be a downstep on the terminal descent whose matching upstep is the middle 
$U$ of a $DUU$.

\section{The bijection} \label{three}

\vspace*{-3mm}

Suppose $(u_i)_{i=1}^{n}$ is a 021-avoiding ascent sequence. Start with $UD$ 
as the current path. For $i=2,3,\dots,n$ in turn, successively increment by 1 the size of the current path $P$ as follows.

\textbf{Case 1.} If $u_i=0$, insert $UD$ at the last peak vertex of $P$, so that 
the resulting path has a long last ascent. 
This is the only case that results in a long last ascent.

\textbf{Henceforth, suppose $\mbf{u_i \ne 0}$.}

\textbf{Case 2.} If $u_i=u_{i-1}$, elevate $P$ (prepend $U$ and append $D$).

Let $a$ and $m$ denote respectively the number of ascents and the 
maximum entry in $(u_1,u_2,\dots,u_{i-1})$.

\textbf{Case 3.} If $u_i=a+1$, append $UD$ to $P$.

\textbf{Case 4.} If we're not in one of the three previous cases, then $u_i \in A_i$ 
(list of allowable $u_i$'s), 
where $A_i :=(m,m+1,\dots,a)$ if $u_{i-1}=0$, and 
$\, :=(m+1,\dots,a)$ if $u_{i-1}> 0$ (which implies that $u_{i-1}=m$, and 
so $u_i=m$ was covered in Case 2). This assertion about $A_i$ holds because the nonzero entries of $(u_1,u_2,\dots,u_{i})$ are weakly increasing and $u_i$ is bounded above by $a+1$. Let $j$ denote the position of $u_i$ in the list $A_i$, 
and $e \ge 0$ the degree of elevation of $P$; 
$e$ is defined because $P$ will not be a pyramid path. 
Insert $UD$ at the top vertex of 
the $j$th key downstep $D_j$ of $P$ and transfer $e$ upsteps from the start of the path to 
the ascent containing the matching upstep of $D_j$. 

The specified insertion is always possible (and reversible) because $|A_i|$ is 
always equal to the  number of key downsteps in $P$, as can be verified by 
a straightforward induction considering the various cases. 
Note that the resulting path has a short last ascent, is not elevated, and does not end with $UD$. \qed

As an example, the 021-avoiding ascent sequence 01012203 produces the following sequence of Dyck paths (key downsteps encountered en route are in blue).
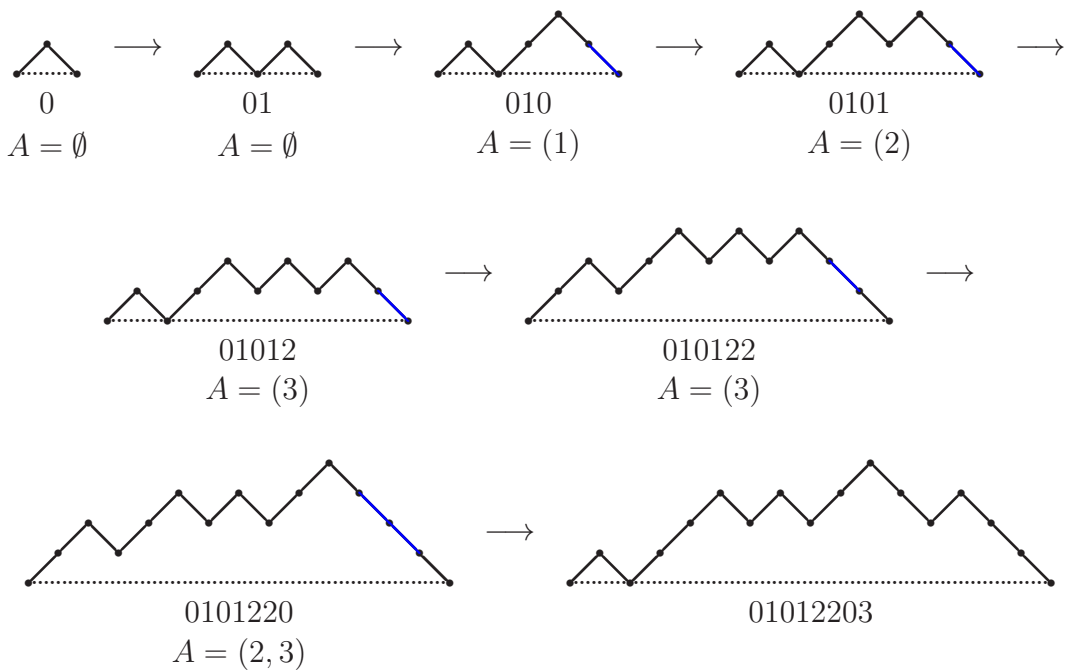
\Einheit=0.4cm
\[
\hspace*{5mm}
\Pfad(-18,0),34\endPfad
\Pfad(-12,0),3434\endPfad
\Pfad(-4,0),343344\endPfad
\Pfad(6,0),34334344\endPfad
\SPfad(-18,0),11\endSPfad
\SPfad(-12,0),1111\endSPfad
\SPfad(-4,0),111111\endSPfad
\SPfad(6,0),11111111\endSPfad
\DuennPunkt(-18,0)
\DuennPunkt(-17,1)
\DuennPunkt(-16,0)
\DuennPunkt(-12,0)
\DuennPunkt(-11,1)
\DuennPunkt(-10,0)
\DuennPunkt(-9,1)
\DuennPunkt(-8,0)
\DuennPunkt(-4,0)
\DuennPunkt(-3,1)
\DuennPunkt(-2,0)
\DuennPunkt(-1,1)
\DuennPunkt(0,2)
\DuennPunkt(1,1)
\DuennPunkt(2,0)
\DuennPunkt(6,0)
\DuennPunkt(7,1)
\DuennPunkt(8,0)
\DuennPunkt(9,1)
\DuennPunkt(10,2)
\DuennPunkt(11,1)
\DuennPunkt(12,2)
\DuennPunkt(13,1)
\DuennPunkt(14,0)
\Label\o{\longrightarrow}(-14,0.2)
\Label\o{\longrightarrow}(-6,0.2)
\Label\o{\longrightarrow}(4,0.2)
\Label\o{\longrightarrow}(16,0.2)
\Label\u{0}(-17,-0.2)
\Label\u{A=\emptyset}(-17,-1.6)
\Label\u{01}(-10,-0.2)
\Label\u{A=\emptyset}(-10,-1.6)
\Label\u{010}(-1,-0.2)
\Label\u{A=(1)}(-1,-1.6)
\Label\u{0101}(10,-0.2)
\Label\u{A=(2)}(10,-1.6)
\blue{ 
\Pfad(1,1),4\endPfad 
\Pfad(13,1),4\endPfad}
\]

\vspace*{3mm}

\Einheit=0.4cm
\[
\hspace*{5mm}
\Pfad(-15,0),3433434344\endPfad
\Pfad(-1,0),334334343444\endPfad
\SPfad(-15,0),1111111111\endSPfad
\SPfad(-1,0),111111111111\endSPfad
\DuennPunkt(-15,0)
\DuennPunkt(-14,1)
\DuennPunkt(-13,0)
\DuennPunkt(-12,1)
\DuennPunkt(-11,2)
\DuennPunkt(-10,1)
\DuennPunkt(-9,2)
\DuennPunkt(-8,1)
\DuennPunkt(-7,2)
\DuennPunkt(-6,1)
\DuennPunkt(-5,0)
\DuennPunkt(-1,0)
\DuennPunkt(0,1)
\DuennPunkt(1,2)
\DuennPunkt(2,1)
\DuennPunkt(3,2)
\DuennPunkt(4,3)
\DuennPunkt(5,2)
\DuennPunkt(6,3)
\DuennPunkt(7,2)
\DuennPunkt(8,3)
\DuennPunkt(9,2)
\DuennPunkt(10,1)
\DuennPunkt(11,0)
\Label\o{\longrightarrow}(-3,0.9)
\Label\o{\longrightarrow}(13,0.9)
\Label\u{01012}(-10,-0.2)
\Label\u{A=(3)}(-10,-1.6)
\Label\u{010122}(5,-0.2)
\Label\u{A=(3)}(5,-1.6)
\blue{ 
\Pfad(-6,1),4\endPfad 
\Pfad(9,2),4\endPfad}
\] 
 
\vspace*{3mm}
 
\Einheit=0.4cm
\[
\Pfad(-17,0),33433434334444\endPfad
\SPfad(-17,0),11111111111111\endSPfad
\Pfad(1,0),3433343433443444\endPfad
\SPfad(1,0),1111111111111111\endSPfad
\DuennPunkt(-17,0)
\DuennPunkt(-16,1)
\DuennPunkt(-15,2)
\DuennPunkt(-14,1)
\DuennPunkt(-13,2)
\DuennPunkt(-12,3)
\DuennPunkt(-11,2)
\DuennPunkt(-10,3)
\DuennPunkt(-9,2)
\DuennPunkt(-8,3)
\DuennPunkt(-7,4)
\DuennPunkt(-6,3)
\DuennPunkt(-5,2)
\DuennPunkt(-4,1)
\DuennPunkt(-3,0)
\DuennPunkt(1,0)
\DuennPunkt(2,1)
\DuennPunkt(3,0)
\DuennPunkt(4,1)
\DuennPunkt(5,2)
\DuennPunkt(6,3)
\DuennPunkt(7,2)
\DuennPunkt(8,3)
\DuennPunkt(9,2)
\DuennPunkt(10,3)
\DuennPunkt(11,4)
\DuennPunkt(12,3)
\DuennPunkt(13,2)
\DuennPunkt(14,3)
\DuennPunkt(15,2)
\DuennPunkt(16,1)
\DuennPunkt(17,0)
\Label\o{\longrightarrow}(-1,1)
\Label\u{0101220}(-10,-0.2)
\Label\u{A=(2,3)}(-10,-1.6)
\Label\u{01012203}(9,-0.2)
\blue{ 
\Pfad(-6,3),44\endPfad 
}
\] 

\vspace*{4mm}
 
The inverse mapping is given in the next section.

\section{The inverse bijection} \label{four}

\vspace*{-3mm}

To reverse the mapping proceed as follows. Start with a ``current path'' taken as 
the given Dyck path. Each step of the algorithm produces an entry of the 
ascent sequence and modifies the current path $P$ to a one-size-smaller 
path according to which of the following four mutually exclusive cases $P$ lies in 
(which match the four cases in Section \ref{three}). Proceed until $P=UD$ and then set $u_1=0$.

In all cases, $i$ denotes the size of the current path $P$, and $Q$ denotes the new 
one-size-smaller path that replaces $P$. An example accompanies the description in each case.

\textbf{Case 1.} The last ascent of $P$ is long. Set $u_{i}=0$. Then delete the last peak.
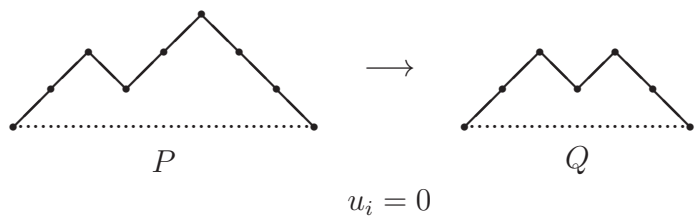
\Einheit=0.5cm
\[
\Pfad(-9,0),33433444\endPfad
\SPfad(-9,0),11111111\endSPfad
\Pfad(3,0),334344\endPfad
\SPfad(3,0),111111\endSPfad
\DuennPunkt(-9,0)
\DuennPunkt(-8,1)
\DuennPunkt(-7,2)
\DuennPunkt(-6,1)
\DuennPunkt(-5,2)
\DuennPunkt(-4,3)
\DuennPunkt(-3,2)
\DuennPunkt(-2,1)
\DuennPunkt(-1,0)
\DuennPunkt(3,0)
\DuennPunkt(4,1)
\DuennPunkt(5,2)
\DuennPunkt(6,1)
\DuennPunkt(7,2)
\DuennPunkt(8,1)
\DuennPunkt(9,0)
\Label\o{\longrightarrow}(1,1)
\Label\u{P}(-5,-0.3)
\Label\u{Q}(6,-0.3)
\Label\u{u_{i}=0}(1,-1.4)
\] 

\vspace*{2mm}

\textbf{Case 2.} The last ascent of $P$ is short and $P$ is elevated. 
Set $u_{i}=u_{i-1}$. (Thus the actual determination of $u_{i}$ is delayed to a 
later step in the algorithm.) Then lower the path, that is, delete the first and last steps.
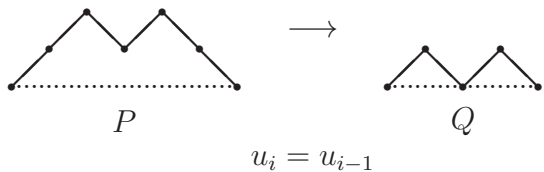
\Einheit=0.5cm
\[
\Pfad(-7,0),334344\endPfad
\SPfad(-7,0),111111\endSPfad
\Pfad(3,0),3434\endPfad
\SPfad(3,0),1111\endSPfad
\DuennPunkt(-7,0)
\DuennPunkt(-6,1)
\DuennPunkt(-5,2)
\DuennPunkt(-4,1)
\DuennPunkt(-3,2)
\DuennPunkt(-2,1)
\DuennPunkt(-1,0)
\DuennPunkt(3,0)
\DuennPunkt(4,1)
\DuennPunkt(5,0)
\DuennPunkt(6,1)
\DuennPunkt(7,0)
\Label\o{\longrightarrow}(1,1)
\Label\u{P}(-4,-0.3)
\Label\u{Q}(5,-0.3)
\Label\u{u_{i}=u_{i-1}}(1,-1.4)
\] 

\vspace*{2mm}

\textbf{Case 3.} The last ascent of $P$ is short and $P$ ends with $UD$. 
(This case is distinct from Case 2 because $P$ has size $\ge 2$.) 
Set $u(i)=$ number of valleys in $P$. Then delete the last peak (= last two steps).

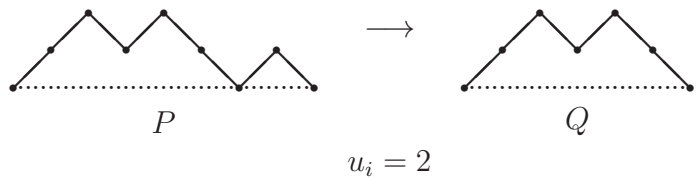
\Einheit=0.5cm
\[
\Pfad(-9,0),33434434\endPfad
\SPfad(-9,0),11111111\endSPfad
\Pfad(3,0),334344\endPfad
\SPfad(3,0),111111\endSPfad
\DuennPunkt(-9,0)
\DuennPunkt(-8,1)
\DuennPunkt(-7,2)
\DuennPunkt(-6,1)
\DuennPunkt(-5,2)
\DuennPunkt(-4,1)
\DuennPunkt(-3,0)
\DuennPunkt(-2,1)
\DuennPunkt(-1,0)
\DuennPunkt(3,0)
\DuennPunkt(4,1)
\DuennPunkt(5,2)
\DuennPunkt(6,1)
\DuennPunkt(7,2)
\DuennPunkt(8,1)
\DuennPunkt(9,0)
\Label\o{\longrightarrow}(1,1)
\Label\u{P}(-5,-0.3)
\Label\u{Q}(6,-0.3)
\Label\u{u_{i}=2}(1,-1.4)
\] 

\vspace*{2mm}

\textbf{Case 4.} The last ascent of $P$ is short and $P$ is neither elevated nor ends $UD$. 
Here we use both the current path $P$ and 
its successor path $Q$ to determine $u_{i}$. Mark the \emph{second} downstep on 
the terminal descent of $P$ and then delete the last peak. The marked downstep 
remains. We need to ensure that it is a key downstep in $Q$. To do so, locate 
its matching upstep, then transfer all upsteps \emph{preceding} this matching 
upstep in its ascent to the start of the path to get $Q$. Note that the marked 
downstep is now indeed a key downstep in $Q$. Set $u_{i} = $ number of valleys 
in $P$ minus the position of the marked downstep among all key downsteps of 
$Q$ when scanned from right to left. This step works because, under the 
bijection, \# ascents in the sequence equals \# valleys in the path 
(proof by induction), and, as noted in Case 4 of the bijection, 
the allowable range $|A_i|$ is 
always equal to the  number of key downsteps in $P$.

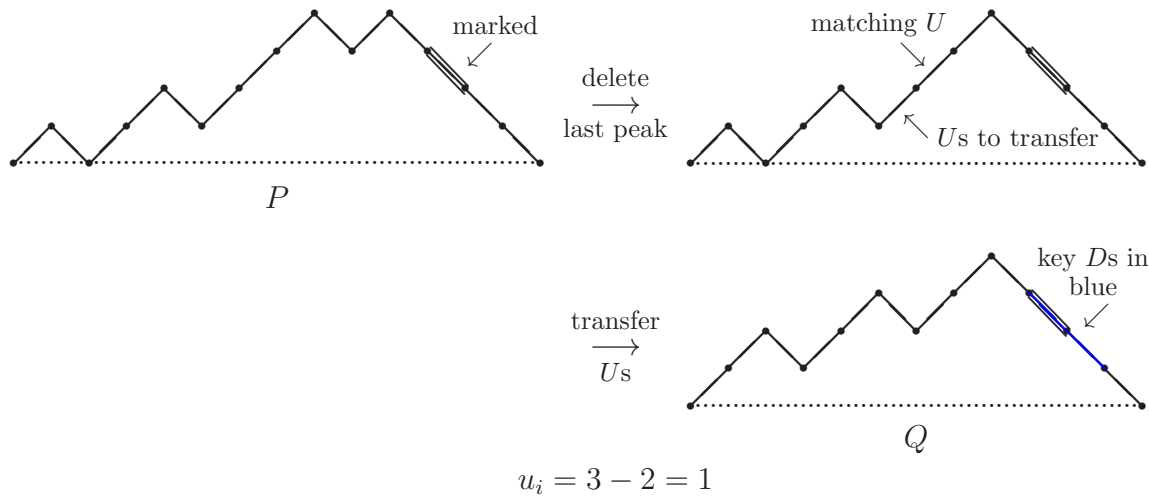
\Einheit=0.5cm
\[
\Pfad(-15,0),34334333434444\endPfad
\SPfad(-15,0),11111111111111\endSPfad
\DuennPunkt(-15,0)
\DuennPunkt(-14,1)
\DuennPunkt(-13,0)
\DuennPunkt(-12,1)
\DuennPunkt(-11,2)
\DuennPunkt(-10,1)
\DuennPunkt(-9,2)
\DuennPunkt(-8,3)
\DuennPunkt(-7,4)
\DuennPunkt(-6,3)
\DuennPunkt(-5,4)
\DuennPunkt(-4,3)
\DuennPunkt(-3,2)
\DuennPunkt(-2,1)
\DuennPunkt(-1,0)
\Pfad(3,0),343343334444\endPfad
\SPfad(3,0),111111111111\endSPfad 
\DuennPunkt(3,0)
\DuennPunkt(4,1)
\DuennPunkt(5,0)
\DuennPunkt(6,1)
\DuennPunkt(7,2)
\DuennPunkt(8,1)
\DuennPunkt(9,2)
\DuennPunkt(10,3)
\DuennPunkt(11,4)
\DuennPunkt(12,3)
\DuennPunkt(13,2)
\DuennPunkt(14,1)
\DuennPunkt(15,0)
\Label\o{\textrm{\footnotesize{delete}}}(1,1.8)
\Label\o{\longrightarrow}(1,1)
\Label\o{\textrm{\footnotesize{last peak}}}(1,0.4)
\Label\u{P}(-8,-0.3)
\Label\u{\textrm{\footnotesize{marked}}}(-2.2,4.3)
\Label\u{\textrm{\footnotesize{$\swarrow$}}}(-2.6,3.5)
\Label\u{\textrm{\Large{$\diagdown$}}}(-3.4,3.2)
\Label\u{\textrm{\Large{$\diagdown$}}}(-3.5,3)
\Label\u{\textrm{\Large{$\diagdown$}}}(12.6,3.2)
\Label\u{\textrm{\Large{$\diagdown$}}}(12.5,3)
\Label\u{\textrm{\footnotesize{$U$s to transfer}}}(11.7,1.3)
\Label\u{\textrm{\footnotesize{$\nwarrow$}}}(9,1.6)
\Label\u{\textrm{\footnotesize{matching $U$}}}(8.0,4.3)
\Label\u{\textrm{\footnotesize{$\searrow$}}}(8.9,3.5)
\]

\Einheit=0.5cm
\[
\white{\Pfad(-15,0),34334333434444\endPfad
\SPfad(-15,0),11111111111111\endSPfad
\DuennPunkt(-15,0)
\DuennPunkt(-14,1)
\DuennPunkt(-13,0)
\DuennPunkt(-12,1)
\DuennPunkt(-11,2)
\DuennPunkt(-10,1)
\DuennPunkt(-9,2)
\DuennPunkt(-8,3)
\DuennPunkt(-7,4)
\DuennPunkt(-6,3)
\DuennPunkt(-5,4)
\DuennPunkt(-4,3)
\DuennPunkt(-3,2)
\DuennPunkt(-2,1)
\DuennPunkt(-1,0)}
\Pfad(3,0),334334334444\endPfad
\SPfad(3,0),111111111111\endSPfad 
\DuennPunkt(3,0)
\DuennPunkt(4,1)
\DuennPunkt(5,2)
\DuennPunkt(6,1)
\DuennPunkt(7,2)
\DuennPunkt(8,3)
\DuennPunkt(9,2)
\DuennPunkt(10,3)
\DuennPunkt(11,4)
\DuennPunkt(12,3)
\DuennPunkt(13,2)
\DuennPunkt(14,1)
\DuennPunkt(15,0)
\Label\u{\textrm{\Large{$\diagdown$}}}(12.6,3.2)
\Label\u{\textrm{\Large{$\diagdown$}}}(12.5,3)
\Label\u{\textrm{\footnotesize{key $D$s in}}}(13.7,4.5)
\Label\u{\textrm{\footnotesize{blue}}}(13.7,3.8)
\Label\u{\textrm{\footnotesize{$\swarrow$}}}(13.8,3.0)
\blue{
\Pfad(12,3),44\endPfad }
\Label\o{\textrm{\footnotesize{transfer}}}(1,1.8)
\Label\o{\longrightarrow}(1,1)
\Label\o{\textrm{\footnotesize{$U$s}}}(1,0.4)
\Label\u{Q}(9,-0.3)
\Label\u{u_{i}=3-2=1}(1,-1.4)
\]
\vspace*{2mm}

\section{Equidistributions} \label{five}

\vspace*{-3mm}

The bijection of Sections \ref{three} and \ref{four} preserves several statistics 
as tabulated below. In the case of the all-zero ascent sequence of length $n$ (which corresponds to the pyramid path $U^n D^n$), a little hiccup arises and the number of terminal 0s must be interpreted as $n-1$.
\begin{center}
  \begin{tabular}{@{} ccc @{}}
    \hline
    021-avoiding ascent sequence 	 & $\longleftrightarrow$ & Dyck path \\[1mm] 
    \hline
    \# initial 0s & $\leftrightarrow$ & length first descent \\ 
    \# terminal 0s & $\leftrightarrow$ & length last ascent $-1$ \\ 
    \# ascents & $\leftrightarrow$ &  $\#\ DU$s   (valleys) \\ 
    \# descents & $\leftrightarrow$ & \# $DUU$s \\ 
    \# entries immediately preceding last    & \raisebox{-1.5ex}{$\leftrightarrow$} & \raisebox{-1.5ex}{degree of elevation}  \\[-3mm] 
    nonzero entry and equal to it & & \\
    \hline
  \end{tabular}
\end{center}

\end{document}